\documentclass[11pt,a4paper]{article}
\usepackage{fullpage}

\usepackage[utf8x]{inputenc} 
\usepackage{amsfonts} 
\usepackage{amsmath,amsthm,amssymb} 

\usepackage{cite} 
\usepackage[numbers]{natbib}

\usepackage{graphicx} 

\usepackage{hyperref}

\usepackage{color}

\newcommand\eps{\varepsilon}
\newcommand\R{\mathbb{R}}
\newcommand\sh{\text{sh}}

\title{Numerical stability of plasma sheath}

\author{Mehdi Badsi$^{1}$, Michel Mehrenberger$^{2,3}$ and Laurent Navoret$^{2,3}$}
\date{}

\begin{document}

\maketitle
\centerline{1. Institut de Math\'ematiques de Toulouse, UMR 5219, }
\centerline{ Universit\'e Toulouse 3, 118 Route de Narbonne,}
\centerline{ 31400 Toulouse, France.}
\bigskip

\centerline{2. Institut de Recherche Math\'ematique Avanc\'ee, UMR 7501,}
\centerline{ Universit\'e de Strasbourg et CNRS, 7 rue Ren\'e Descartes,}
\centerline{ 67000 Strasbourg, France}
\bigskip

\centerline{3. INRIA Nancy-Grand Est, TONUS Project, Strasbourg, France}

\bigskip

%
%
\begin{abstract}We are interested in developing a numerical method for capturing stationary sheaths, that a plasma forms in contact with a metallic wall. This work is based on a bi-species (ion/electron) Vlasov-Amp\`ere model proposed in \cite{2016_Badsi}. The main question addressed in this work is to know if classical numerical schemes can preserve stationary solutions with boundary conditions, since these solutions are not a priori conserved at the discrete level.  In the context of high-order semi-Lagrangian method, due to their large stencil, interpolation near the boundary of the domain requires also a specific treatment.  \end{abstract}

\maketitle
\section{Introduction}

When a plasma is in contact with a metallic wall, stationary boundary layers, called sheaths, form. They are due to both the reflexion/absorption properties of the wall and the charge imbalance coming from the mass difference between electrons and ions: electrons leave the plasma faster than ions. A drop of the self-consistent potential at the wall is indeed built up to accelerate ions and decelerate electrons just as to ensure an equal flux of ions and electrons (a zero current) at the wall. In this work, the starting point is a plasma sheath model proposed by Badsi et al. \cite{2016_Badsi}, for which a so-called kinetic Bohm criterion on the incoming ion distribution has been found out. Numerically, the sheath potential is computed by a non-linear Poisson solver. 

We here consider the corresponding \emph{non-stationary} ion/electron Vlasov-Amp\`ere model (note that the stability of its linearized version has been studied in \cite{BADSI2017954}) and its  numerical implementation through a semi-Lagrangian scheme \cite{Sonnen}. 
Semi-Lagrangian schemes are transport solvers with no stability conditions on the time step: time step is only constrained by the physical dynamics that has to be captured. High-order schemes, using high degree interpolation and time splitting, can be devised  and parallel implementation can be considered.  We raise here the following questions:
\begin{itemize}
\item Question 1: How do semi-Lagrangian scheme preserve this equilibrium in time ? 
\item Question 2: How can we treat boundary conditions with large stencil interpolation ?
\end{itemize}

We will see that the equilibrium is stable provided that the incoming current, not exactly zero at the discrete level,  is made small. 
This will be ensured by taking a fine grid in space and velocity and a time step small enough. 
For the second question, fictitious values of the distributions functions will be set at the spatial boundaries.
Finally, we emphasize that the time step is constrained by the electronic dynamics and thus is very small here.

\section{Sheath model}

\subsection{Stationary solution}

 Following \cite{2016_Badsi}, we consider a one-dimensional domain $[0, 1]$: particles enter the domain on the left ($x = 0$) and the metallic wall is located at abscissa $x=1$. We consider in this work purely absorbing boundary conditions at the wall.

Electrons and ions are described by their distribution functions in phase space, denoted $f_{s_{e}}^{\sh} (x,v)$ and $f_{s_{i}}^{\sh} (x,v)$ with spatial variable $x \in [0,1]$ and velocity variable $v \in \R$, satisfy the dimensionless stationary Vlasov system:
\begin{align}
&v\, \partial_{x} f_{s_{e}}^{\sh}  - \frac{1}{\mu} E^{\sh} \, \partial_{v} f_{s_{e}}^{\sh} = 0,\label{eq:fe_stationary}\\ 
&v\, \partial_{x} f_{s_{i}}^{\sh}  + \phantom{\frac{1}{\mu}} E^{\sh}\, \partial_{v} f_{s_{i}}^{\sh}  = 0,\label{eq:fi_stationary}
\end{align}
where $\mu \ll 1$ is the mass ratio between electrons and ions, and $E^{\sh}(x)$ denotes the stationary electric field. It satisfies the Poisson equation:
\begin{align}
&E^{\sh} = -\frac{d \phi}{dx}^{\sh},\quad  -\eps^{2}\, \frac{d^{2} \phi}{dx^{2}}^{\sh} = \rho^{\sh},\label{eq:Poisson}\\
&\rho^{\sh}(x) = \int_{v\in \R} (f_{s_{i}}^{\sh} (x,v)- f_{s_{e}}^{\sh} (x,v))\, dv,
\end{align}
where $\phi^{\sh}(x)$ is the electric potential, $\eps \ll 1$ is the dimensionless Debye length of the plasma and $\rho^{\sh}(x)$ is the charge density.

\noindent {\bfseries Boundary conditions.} At the entrance of the domain, incoming ions and electrons have distributions:
\begin{equation}
\forall v > 0,\quad f_{s_{i}}^{\sh}(0,v) = f_{s_{i}}^{in}(v), \quad \quad f_{s_{e}}^{\sh}(0,v) = n_{0}\, \sqrt{\frac{2\mu}{\pi}} e^{-\displaystyle\frac{\mu v^{2}}{2}},\label{eq:initdistrib}
\end{equation}
where $n_{0}$ is an electronic density parameter, that will be determined to have a given charge density $\rho_{0}$ at the entrance of the domain. Thus, electrons have a Maxwellian distribution. At the wall, ions and electrons are supposed here to be totally absorbed:
\begin{equation*}
\forall v < 0,\quad f_{s_{i}}^{\sh}(1,v) = 0, \quad \quad f_{s_{e}}^{\sh}(1,v) = 0,
\end{equation*}
The boundary conditions of the potential \eqref{eq:Poisson} are defined as follows. The potential values at $x=0$ is set to $0$ and the potential at the wall, 
$\phi_{w} = \phi^{\sh}(1)$, also called floating potential, is fixed such that the current vanishes at the wall:
\begin{align}
&\phi^{\sh}(0) = 0,\label{eq:Poisson_dirichlet}\\
&\int_{v\in \R} v f_{s_{i}}^{in}(v)\, dv = \sqrt{\frac{2}{\pi\mu}}n_{0}e^{\phi_{w}}.\label{eq:zero_currrent}
\end{align} 
 Note that it implies that the current vanishes in the whole domain, since both electron and ion momentum are constant in space.

\noindent {\bfseries Solutions.}  Using the characteristic lines of the transport equation and under the hypothesis that the potential is decreasing in the domain, the solution to problem \eqref{eq:fe_stationary}-\eqref{eq:fi_stationary}-\eqref{eq:Poisson} is given by:
\begin{align}
f_{s_{i}}^{\sh}(x,v) = &
\boldsymbol{1}_{\left\{v> \sqrt{-2\,\phi^{\text{sh}}(x)} \right\}}\, f_{s_{i}}^{in}\Big(\sqrt{v^{2} + 2\,\phi^{\text{sh}}(x)}\Big),\label{eq:ion_sol}\\
f_{s_{e}}^{\text{sh}}(x,v) = &
\boldsymbol{1}_{\left\{v \geq -\sqrt{\frac{2}{\mu}\,(\phi^{\text{sh}}(x)-\phi^{\text{sh}}(1))}\right\}}\, f_{s_{e}}^{in}\Big(\sqrt{v^{2} - \frac{2}{\mu}\phi^{\text{sh}}(x)}\Big),\label{eq:elec_sol}
\end{align}
where the sheath potential is the solution to the following non-linear Poisson equation:
\begin{equation}
-\eps^{2}\, \frac{d^{2} \phi}{dx^{2}}^{\sh}(x) = \left[\int_{\R^{+}} \frac{f_{s_{i}}^{in}(v) v}{\sqrt{v^{2} - 2\,\phi^{\sh}(x)}}\, dv\right] - \frac{2n_{0}}{\sqrt{2\pi}}\left[\sqrt{2\pi}\, e^{\phi^{\sh}(x)} -  \int_{\sqrt{-2\,\phi_{w}}}^{+\infty} \frac{e^{-\frac{ v^{2}}{2}} v}{\sqrt{v^{2} + 2\,\phi^{\sh}(x)}}\, dv\right].\label{eq:Poisson_non_linear}
\end{equation} 
complemented by the boundary conditions \eqref{eq:Poisson_dirichlet}-\eqref{eq:zero_currrent}.

Let consider a given charge density $\rho_{0} \in \R$ at the entrance of the domain.  The problem is solved in two steps:
\begin{enumerate}
\item a negative potential drop at the wall $\phi_{w} \leq 0$ is uniquely defined satisfying \eqref{eq:zero_currrent} provided that the ionic distribution satisfies the relation:
\begin{equation*}
\frac{\int_{\R_{+}} f_{s_{i}}^{in}(v) v\, dv}{\int_{\R_{+}} f_{s_{i}}^{in}(v)\, dv - \rho_{0}} \leq \sqrt{\frac{2}{\mu\pi}}.
\end{equation*}
It satisfies the following equation:
\begin{equation}
\frac{1}{\sqrt{\mu}} e^{\phi_{w}} \left(\int_{v\in \R} f_{s_{i}}^{in}(v)\, dv - \rho_{0}\right)+ \left(\int_{v\in \R} v f_{s_{i}}^{in}(v)\, dv\right) \left(\int_{\sqrt{-2\phi_{w}}}^{+\infty} e^{-\frac{v^{2}}{2}}\, dv\right) = 0,\label{eq:phiw}
\end{equation}
and the associated electronic density parameter is given by:
\begin{equation*}
n_{0} = \sqrt{\frac{\pi}{2}} \frac{\int_{v\in \R_{+}} f_{s_{i}}^{in}(v)\, dv - \rho_{0}}{(\sqrt{2\pi} - \int_{\sqrt{-2\phi_{w}}}^{+\infty} e^{-\frac{v^{2}}{2}}\, dv)}.\label{eq:n0}
\end{equation*}
\item then, once this potential wall exists, the non-linear Poisson equation \eqref{eq:Poisson_non_linear}, complemented with the boundary conditions $\phi^{\sh}(0) = 0$ and $\phi^{\sh}(1) = \phi_{w}$, has a unique solution under the kinetic Bohm criterion: 
\begin{equation*}
\frac{\displaystyle\int_{\R_{+}} \frac{f_{s_{i}}^{in}(v)}{v^{2}}\,dv}{\displaystyle\int_{\R_{+}} f_{s_{i}}^{in}(v)\,dv} < \frac{\sqrt{2\pi} + \displaystyle\int_{\sqrt{-2\phi_{w}}}^{+\infty} \frac{e^{-v^{2}/2}}{v^{2}}\,dv}{\sqrt{2\pi} - \displaystyle\int_{\sqrt{-2\phi_{w}}}^{+\infty} e^{-v^{2}/2}\,dv}.
\end{equation*}
\end{enumerate}

\noindent {\bfseries Parameters.} The equilibria has been computed by the gradient algorithm described in \cite{2016_Badsi} with the following parameters: 
\begin{itemize}
\item the dimensionless Debye length is set to $\eps = 0.01$.
\item the charge is zero at the entrance of the domain: $\rho_{0} = 0$.
\item the incoming ion distribution is given by
\begin{equation}
f_{s_{i}}^{in}(v) = \boldsymbol{1}_{\left\{v> 0\right\}}\,\min(1,v^{2}/\eta) \frac{1}{\sqrt{2\pi \sigma^{2}}} e^{- \frac{(v-Z)^{2}}{2\sigma^{2}}},\qquad \eta = 10^{-1},\quad \sigma = \sqrt{T_{i}/T_{e}} = 0.5,\quad Z=1.5.\label{eq:ion_income}
\end{equation}
$Z$ is the macroscopic ion velocity when entering the domain.
\item the incoming electronic distribution is Maxwellian as given by \eqref{eq:initdistrib}. For the sake of completeness, we recall the expression:
\begin{equation}
f_{s_{e}}^{in}(v) =  \boldsymbol{1}_{\left\{v> 0\right\}}\, n_{0}\, \sqrt{\frac{2\mu}{\pi}} e^{-\displaystyle\frac{\mu v^{2}}{2}},\qquad \mu = 1/3672.
\end{equation}
$\mu$ is the mass ratio of a deuterium plasma. 
\end{itemize}
Given these parameters, we can compute $n_0,\phi_w$ with \eqref{eq:phiw}-\eqref{eq:n0}, and the sheath electric potential $\phi^{\text{sh}}$ by solving the non-linear Poisson equation \eqref{eq:Poisson_non_linear}. 
In the numerical results, we take: 
$$
n_0 =  0.50191266314760252,\ \phi_w = -2.7839395640524267.
$$
In practice, for future use in the non-stationary model, we store the nonlinear potential $\phi^{\text{sh}}$  on a uniform grid of $N=2048$ cells, that is, at points $j/N,\ j=0,\dots,N$.

\subsection{Non-stationary model}

In the non-stationary setting, electron and ion distribution functions, $f_{s_{e}}(t,x,v)$ and $f_{s_{i}}(t,x,v)$, and the electric field, $E(t,x)$, depend also on time $t > 0$.
They satisfy the non-linear Vlasov-Amp\`ere system:
\begin{align}
&\partial_{t} f_{s_{e}} + v\, \partial_{x} f_{s_{e}}  - \frac{1}{\mu} E\, \partial_{v} f_{s_{e}} = 0,\label{eq:fe}\\ 
&\partial_{t} f_{s_{i}} + v\, \partial_{x} f_{s_{i}} + \phantom{\frac{1}{\mu}}  E\, \partial_{v} f_{s_{i}} = 0,\label{eq:fi}\\
&\eps^{2}\,\partial_{t} E = -J,\label{eq:E}
\end{align}
where $J$ denotes the current density:
\begin{equation*}
J(t,x) =   \int_{v\in \R} v\, (f_{s_{i}} (t,x,v) - f_{s_{e}} (t,x,v))\, dv.
\end{equation*}

\noindent {\bfseries Initial data.} The initial data are given by the stationary solution:\begin{equation*}
f_{s_{e}}(0,x,v)  = f_{s_{e}}^{\sh}(x,v),\quad f_{s_{i}}(0,x,v) = f_{s_{i}}^{\sh}(x,v),\quad E(0,x) = E^{\sh}(x).
\end{equation*}
We note that the initial electric field $E(x,0)$ satisfies the following Poisson equation:
\begin{align}
&\eps^{2}\partial_{x} E(0,x) = n_{s_{i}}(0,x) - n_{s_{e}}(0,x),\label{eq:Einit}\\
&\int_{0}^{1} E(t,x)\, dx = - \phi_{w},\label{eq:Einit2}
\end{align}
where the densities $n_{s_{i}}$ and $n_{s_{e}}$ are given by:
\begin{equation*}
n_{s_{e}}(t,x) = \int_{v\in \R} f_{s_{e}} (t,x,v)\, dv,\quad n_{s_{i}}(t,x) = \int_{v\in \R} f_{s_{i}} (x,v)\, dv,
\end{equation*}
and $\phi_{w}$ is the floating potential at $x = 1$.

\section{Numerical scheme}

We solve system \eqref{eq:fe}-\eqref{eq:fi}-\eqref{eq:E} using a semi-Lagrangian scheme. Due to the large mass ratio between electrons and ions, the computational domain of electrons will be larger than the ion one and thus velocity meshes have to be different for the two kinds of particles. 

\noindent {\bfseries Mesh notations.} We consider uniform cartesian meshes for both spatial domain $[0,1]$ and velocity domain $[v_{\min},v_{\max}]$.
Considering $N_{x}+1$ points in the spatial direction and $N_{v}+1$ points in the velocity one, mesh points 
are denoted $(x_{i}, v_{j}) = (i\Delta x, v_{\min}+j\Delta v)$ 
 for all $0\leq i \leq N_{x}$ and $0\leq j\leq N_{v}$ with $\Delta x = 1/N_{x}$ and 
$\Delta v = (v_{\max}- v_{\min}) / N_{v}$.
Note that the definitions relative to the velocity domain 
implicitly depend on the particle $s\in \{s_e,s_i\}$ under consideration. 
We denote by $f_{s, (i,j)}^{n}$ the approximate value of the distribution function $f_{s}$ at point $(x_{i}, v_{j})$ at time $t_{n} = n\Delta t$, with $\Delta t > 0$ and $n \in \mathbb{N}$, and $E_{i}^n$ the approximate value of the electric field at point $x_{i}$. 

\subsection{Initialization of the electric field}

The solution to system \eqref{eq:Einit}-\eqref{eq:Einit2} is given by:
\begin{equation}
E(x) = \int_{0}^{x}  \frac{n_{s_{i}}(y) - n_{s_{e}}(y)}{\eps^{2}}\, dy  -\phi_{w} - \int_{0}^{1} \int_{0}^{x}  \frac{n_{s_{i}}(y) - n_{s_{e}}(y)}{\eps^{2}}\, dy dx.\label{eq:Esol}
\end{equation}
The charge densities $n_{s}, s\in \{s_e,s_i\}$ are computed at grid points $x_{i}$ using the trapezoidal formula in the velocity direction from the distributions functions \footnote{The trapezoidal formula is spectrally accurate for smooth periodic data and remains
accurate when the velocity domain is large enough so that the distribution equals zero at the boundaries (up to machine precision) and can be considered as periodic. We yet point out that the electron distribution \eqref{eq:elec_sol} has a discontinuity in velocity, which deteriorates the accuracy.}. We then consider a reconstruction at any point $x = x_{i}+\alpha \Delta x \in [0,1]$, with $\alpha \in [0,1]$,  by local centered Lagrange interpolation of degree $2d+1$ with $d \in \mathbb{N}$:
\begin{equation*}
\quad n_{s, h}(x) = \sum_{k=-d}^{d+1} n_{s, i+k}\, L_{k}(\alpha),
\end{equation*}
where $L_{k} $ are the elementary Lagrange polynomials: 
\begin{equation}L_{k}(\alpha) =\prod_{\begin{subarray}{l} i=-d\\
i \neq k
\end{subarray}}^{d+1} (\alpha-i)/(k-i).\label{eq:Lk}
\end{equation} 
Note that the densities have to be defined at points outside of the domain when $d > 0$. This will be done by expanding the definition of the distribution functions. This point will be detailed in Section \ref{bc}. The discrete electric field at grid points $x_{i}$ is then obtained from \eqref{eq:Esol} in which $n_{s}$ are replaced by their discrete counterparts $n_{s,h}$. The involved integrals are computed exactly and the overall scheme has $O(\Delta x^{2d+2})$ accuracy.

\subsection{Splitting}

To solve the Vlasov-Amp\`ere system, we use a splitting between space and velocity dynamics and write it as a succession of one-dimensional advections.  More precisely, we consider the following dynamics: the kinetic transport system given by:
\begin{align}
(\mathcal{T})\qquad &\partial_{t} f_{s_{e}} + v\,\partial_{x} f_{s_{e}} = 0,\label{eq:kine_transp_x}\\
\qquad &\partial_{t} f_{s_{i}} + v\,\partial_{x} f_{s_{i}} = 0,\label{eq:kin_transp_x2}\\
&\eps^{2}\,\partial_{t} E = -J,\label{eq:E_x}
\end{align}
and the electric transport system, given by:
\begin{align}
(\mathcal{U})\qquad&\partial_{t} f_{s_{e}} -\frac{1}{\mu} E\,\partial_{v} f_{s_{e}} = 0,\label{eq:kin_transp_v}\\
&\partial_{t} f_{s_{i}} + E\,\partial_{v} f_{s_{i}} = 0,\label{eq:kin_transp_v2}\\
&\eps^{2}\,\partial_{t} E = 0.\label{eq:E_v}
\end{align}
Since the electric field is constant in time in this second step, both dynamics $\mathcal{T}$ and $\mathcal{U}$ are advections at constant velocities.  
In practice, to obtain second order accuracy in time, we consider the Strang splitting which consists in computing for $s\in \{s_e,s_i\}$ 
\begin{equation*}
\left\{\left(f_{s,(i,j)}^{n+1}\right)_{i,j}, (E_{i}^{n+1})_{i}\right\} = \Big[\mathcal{U}_{h,\Delta t/2}\circ\mathcal{T}_{h,\Delta t}\circ\mathcal{U}_{h,\Delta t/2}\Big] \left\{ (f_{s,(i,j)}^{n})_{i,j}, (E_{i}^{n})_{i}\right\}
\end{equation*}
where $\mathcal{T}_{h,\tau}$ and $\mathcal{U}_{h,\tau}$ are discrete approximations of $\mathcal{T}$ and $\mathcal{U}$ over a time interval $\tau$.

Each transport equation is solved using a semi-Lagrangian scheme with centered Lagrange interpolation of degree $2d+1$. We thus consider the discretized version of the kinetic and electric transport dynamics. Operator $\mathcal{T}_{h,\tau}: ((f_{s, (i,j)}), (E_{i})_{i}) \rightarrow ((f_{s, (i,j)}^{\ast}), (E^{\ast}_{i})_{i})$ consists in the following: for any grid points $(x_{i},v_{j})$, we define the shifted index $i^\ast$ and $\alpha \in [0,1]$ such that $x_{i} - v_{j}\tau = x_{i^{\ast}} + \alpha \Delta x$ and then the distribution function is given by:
\begin{equation*}
f_{s,(i,j)}^{\ast} = \sum_{k=-d}^{d+1} f_{s,(i^{\ast}+k,j)} L_{k}(\alpha),
\end{equation*}
where the $(L_{k})$ are defined in \eqref{eq:Lk}, and the electric field at point $x_{i}$ is given by:
\begin{equation*}
E^{\ast}_{i} = E_{i} - \tau\, \frac{J_{i}+J_{i}^{\ast}}{2},
\end{equation*}
where the current density $(J_{i})_{i}$ (resp. $(J_{i}^{\ast})_{i}$) is computed by trapezoidal rule in velocity from the discrete distribution function $(f_{i,j})$ (resp. $(f_{i,j}^{\ast})$). We thus exactly solve in time the transport equations \eqref{eq:kine_transp_x}-\eqref{eq:kin_transp_x2} at the grid points, starting from the interpolated distribution function in the spatial direction, while the Amp\`ere equation \eqref{eq:E_x} is computed using a second order Crank-Nicolson scheme. The interpolation requires values of the distribution function outside the domain: we explain in the next section how to extrapolate it.

Operator $\mathcal{U}_{h,\tau}: ((f_{s, (i,j)}), (E_{i})_{i}) \rightarrow ((f_{s, (i,j)}^{\ast}), (E^{\ast}_{i})_{i})$ consists in the following: for any grid points $(x_{i},v_{j})$, we define the shifted indexes $j^\ast_{s_{i}}$, $j^{\ast}_{s_{e}}$ and $\alpha_{s_{i}}, \alpha_{s_{e}} \in [0,1]$ such that $v_{j} - E_{i}\tau = v_{j^{\ast}_{s_{i}}} + \alpha_{s_{i}} \Delta v_{s_{i}}$ and $v_{j} + \frac{1}{\mu} E_{i}\tau = v_{j^{\ast}_{s_{e}}} + \alpha_{s_{e}} \Delta v_{s_{e}}$ and then
\begin{align*}
f_{s,(i,j)}^{\ast} &= \sum_{k=-d}^{d+1} f_{s,(i,j^{\ast}_{s}+k)} L_{k}(\alpha_{s}), \\
E^{\ast}_{i} &= E_{i}.
\end{align*}
For this advection in velocity $\mathcal{U}_{h,\tau}$,  periodic boundary conditions are used; we have here made the presentation for Lagrange interpolation, as it is used for advection in space,
but other advection scheme can be used; in particular, in the numerical results, we will use cubic splines.

\subsection{Boundary conditions}
\label{bc}

In both the initial computation of the electric field and the advection in space $\mathcal{T}_{h,\tau}$, the proposed numerical scheme requires to take values of the distribution function outside the physical domain.

For any $x_{i} = i\Delta x < 0$ and $v_{j}$, we consider the following extrapolation at the entry $x=0$:
\begin{equation*}
f_{s, (i,j)} = \begin{cases}
f_{s}(0,0,v_j),&\text{if } v_{j} \geq 0,\\
2 f_{s, (0,j)} - f_{s, (-i,j)} ,&\text{if } v_{j} < 0.\\
\end{cases}
\end{equation*}
For any $x_{i+N_{x}} = (i+N_{x})\Delta x > 1$ and $v_{j}$, we consider the following extrapolation at the wall $x=1$:
\begin{equation*}
 f_{s, (i+N_{x},j)} = \begin{cases}
2 f_{s, (N_{x},j)} - f_{s, (N_{x}-i,j)} &\text{if } v_{j} \geq 0,\\
0&\text{if } v_{j} < 0\\
\end{cases}.
\end{equation*}
This corresponds to a purely Dirichlet condition for incoming velocities and an extension by imparity for leaving velocities (also called \emph{butterfly} procedure).


\section{Numerical results}

\subsection{Sheath test-case}

We here used the following set of parameters: 
\begin{align}
&d = 8,\quad N_{x} = 2048,\quad N_{v} = 4096,\nonumber\\
&\text{velocity domain } [−200, 500] \text{ for electrons and } [-5,5] \text{ for ions},\label{parameters}\\ 
&\Delta t = 10^{-5}.\nonumber
\end{align}
The simulation has run on $256$ processors during $24$ hours, with final time $t=8.03478$, on the Marconi supercomputer; the distribution function is stored every multiple of $0.01$,
and the time diagnostics every time step.

\begin{figure}
\begin{tabular}{c}
\includegraphics[width=1.2\linewidth]{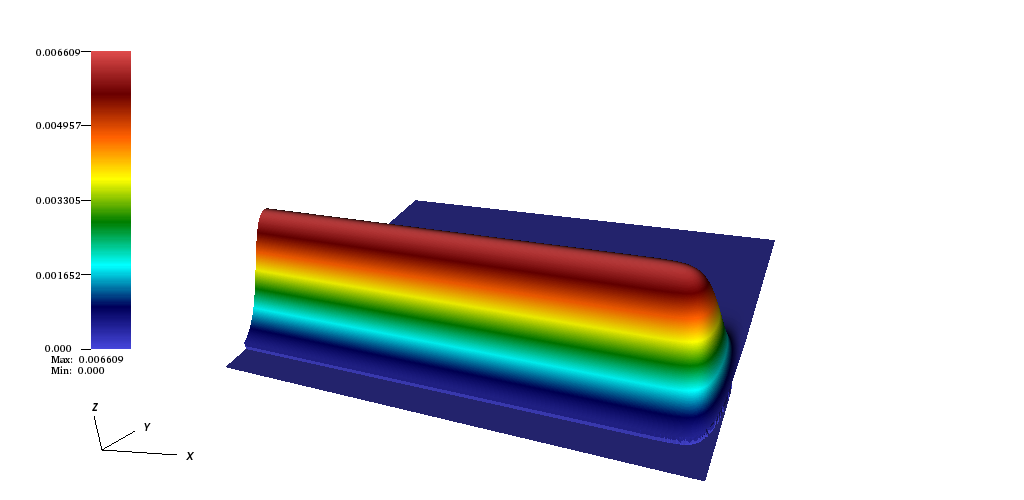} \\
\includegraphics[width=1.2\linewidth]{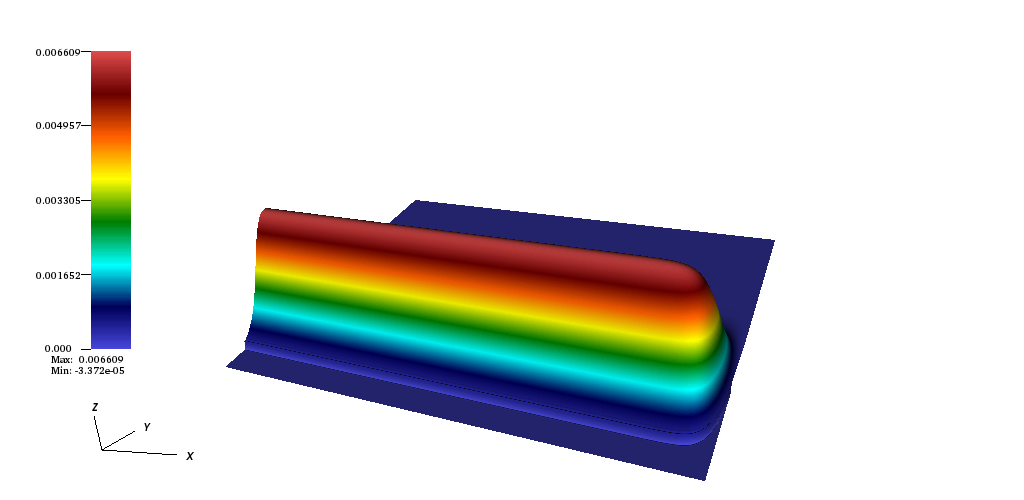}
\end{tabular}
\caption{Distribution function of electrons at time $t=0$ (top) and $t=4$ (bottom);  X stands for $x$, Y for $v$ and Z for $f(t,x,v)$. Parameters given in \eqref{parameters}.}
\label{fig1}
\end{figure}

\begin{figure}
\begin{tabular}{cc}
\includegraphics[width=1.2\linewidth]{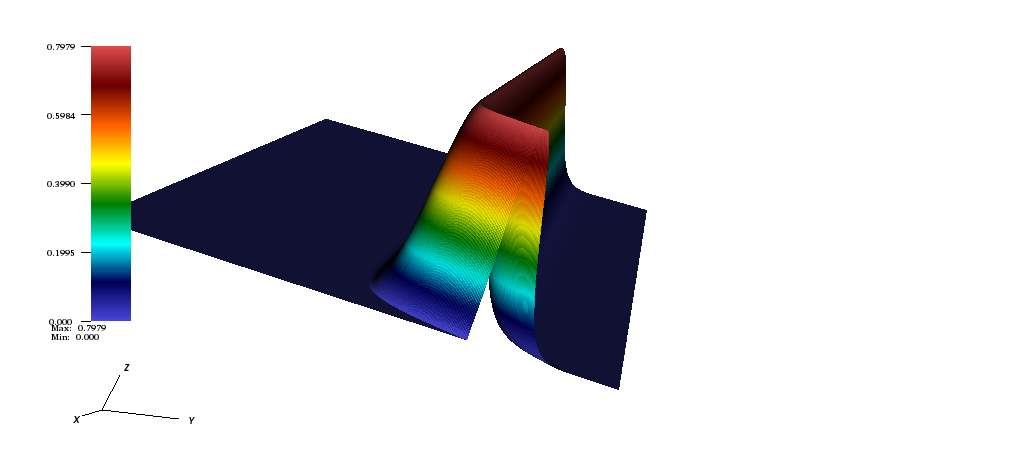} \\
\includegraphics[width=1.2\linewidth]{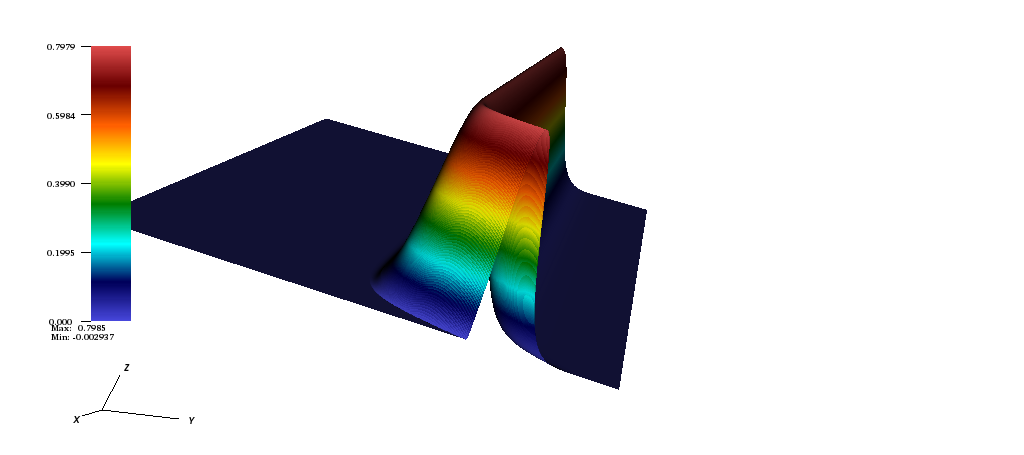}
\end{tabular}
\caption{Distribution function of ions at time $t=0$ (top) and $t=4$ (bottom); X stands for $x$, Y for $v$ and Z for $f(t,x,v)$. Parameters given in \eqref{parameters}.}
\label{fig2}
\end{figure}

\textbf{Distribution functions.} On Figure \ref{fig1} (resp. Figure \ref{fig2}), we represent the distribution of electrons (resp. ions) at time $t=0$ (top) and $t=4$ (bottom).
We clearly see that the equilibrium is well preserved. 
 We see that the maximum principle is not exactly preserved, as no limiting procedure is introduced both in space advection ($d=8$, i.e. Lagrange interpolation of degree $17$)
and in velocity advection (cubic splines). However, this seems to be not crucial here, as it is not far from being preserved and it can be used as a measure of accuracy of the simulation as other theoretically 
preserved quantities. We see that the distribution function of electrons presents a discontinuity; this has the effect that the quadrature in velocity for computing the current converges slowly, and thus a high number of points in velocity is needed.

For ions, we see that the distribution is not constant in space near the wall ($x=1$) for a given velocity. At this boundary, the butterfly procedure does not destroy the $C^1$ property of the distribution function and seems to be more adapted than the prolongation by a constant value, which is used for incoming velocities. Note also that the fictitious boundary values are only used to interpolate inside the domain, as the sub time step $\tau$ is here always positive
(negative time steps, could be however considered when going to higher order splitting schemes in time, but are not studied in this work). 

The space discretization is rather fine; this is needed for the sharp gradient near the wall. Non-uniform grids could be useful here to save memory and computations \cite{CouletteManfredi2014}, 
but are not tackled here for simplicity. We mention also that high order interpolation is used, which permits not to have to refine too much. High order schemes in a non uniform setting could be considered
with a Semi-Lagrangian Discontinuous Galerkin method (SLDG) (see \cite{SLDG-CMV2010}, for example), or with non uniform cubic splines (see for example \cite{Afeyan2014}). 

 Finally, we note that the time step is chosen very small. In order to capture the oscillations due to the dynamic of the electrons, it should be indeed smaller than $\simeq 2\cdot 10^{-4}$.
However, as we are interested in studying the stability of the equilibrium for long times (at least \emph{long} with respect to the electron dynamics), we have to reduce the time step to one order of magnitude
in order to reduce the time error (due to the Strang splitting and the Crank-Nicolson scheme for the Amp\`ere equation). Here again, higher order schemes in time might be useful, in order to use larger time steps
(see for example \cite{Afeyan2014}).

\begin{figure}
\begin{tabular}{cc}
\includegraphics[width=0.4\linewidth]{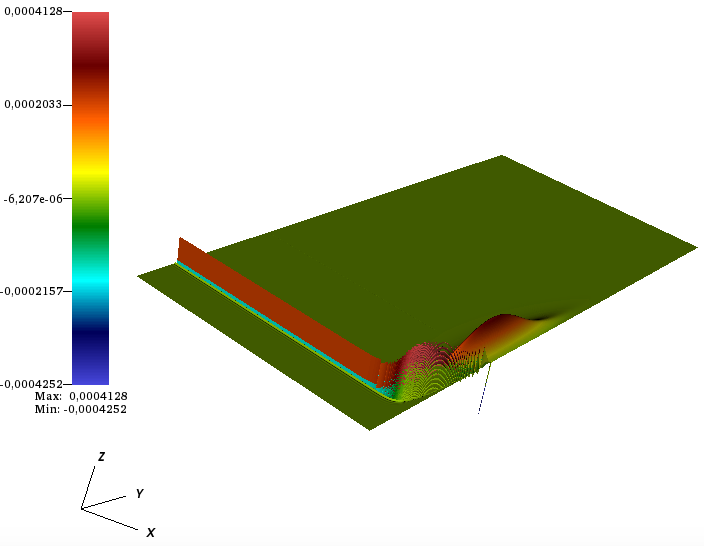} &
\includegraphics[width=0.4\linewidth]{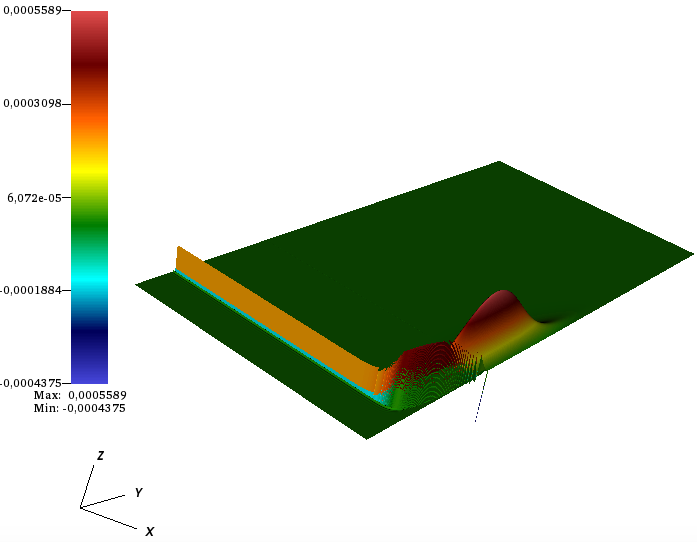} 
\end{tabular}
\caption{Error  on the electron distribution function, $f_{s_e}(0,x,v)-f_{s_e}(t,x,v)$, at time $t=4$ (left) and $t=8$ (right);  X stands for $x$, Y for $v$ and Z for $f_{s_e}(0,x,v)-f_{s_e}(t,x,v)$. Parameters given in \eqref{parameters}.}
\label{fig3}
\end{figure}

\begin{figure}
\begin{tabular}{cc}
\includegraphics[width=0.4\linewidth]{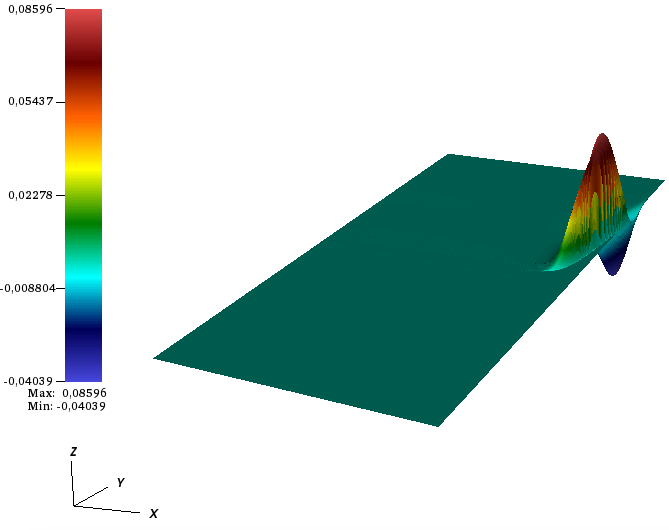} &
\includegraphics[width=0.4\linewidth]{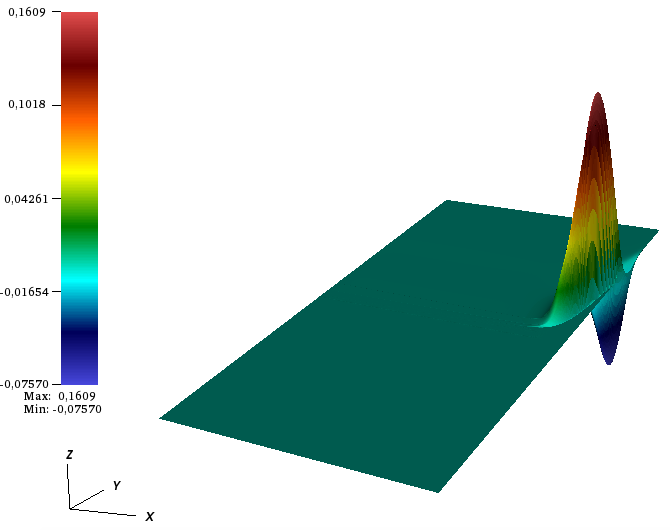} 
\end{tabular}
\caption{Error  on the ion distribution function, $f_{s_i}(0,x,v)-f_{s_i}(t,x,v)$, at time $t=4$ (left) and $t=8$ (right);  X stands for $x$, Y for $v$ and Z for $f_{s_i}(0,x,v)-f_{s_i}(t,x,v)$. Parameters given in \eqref{parameters}.}
\label{fig4}
\end{figure}

On Figure \ref{fig3} (resp. \ref{fig4}), we represent the difference between the distribution function of electrons (resp. ions) at time $0$ and at time $t=4$ and $t=8$. We see that the error for the electrons is mainly localized near the wall and at the discontinuity of the distribution.
As regards ions, the error is localized near the wall. We see that the error near the wall increases in time, especially for the ions. We remark also an exceptional value for the electrons, that seems
to be located at $v=0$ and $x=1$: this might be explained by the fact that the function is not modified by the space advection (at $v=0$), and is thus not diffused.  

\begin{figure}
\begin{tabular}{cc}
\includegraphics[width=0.4\linewidth]{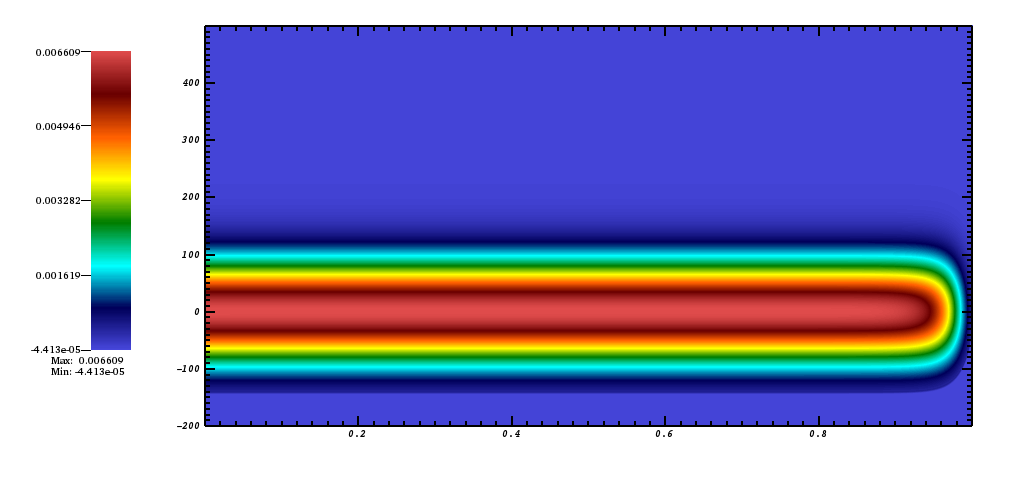} &
\includegraphics[width=0.4\linewidth]{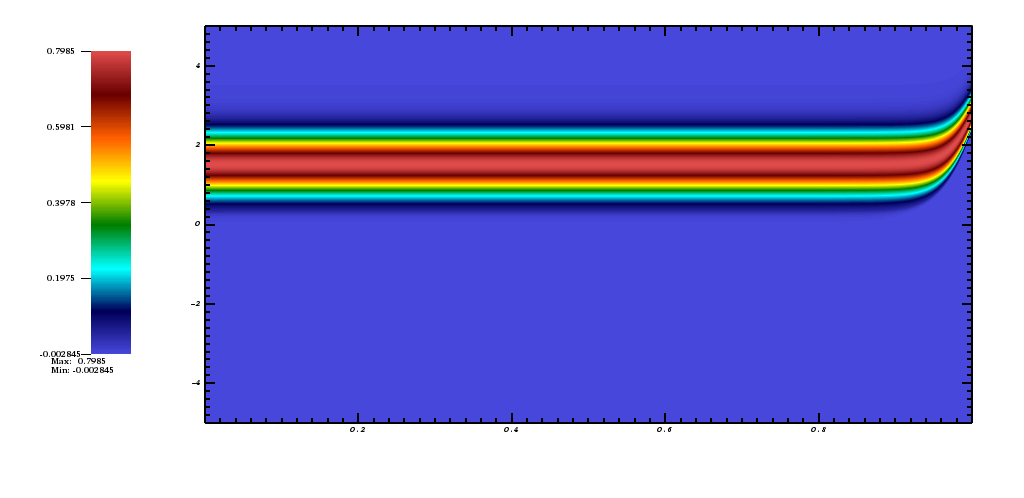} 
\end{tabular}
\caption{Distribution function of electrons (left) and ions (right) at time $t=8$. Parameters given in \eqref{parameters}.}
\label{fig5}
\end{figure}

On Figure \ref{fig5}, we represent $2d$ views of the distribution function of electrons and ions at time $t=8$. We recognize here the pictures of \cite{2016_Badsi}.

\begin{figure}
\centering
\includegraphics[width=0.8\linewidth]{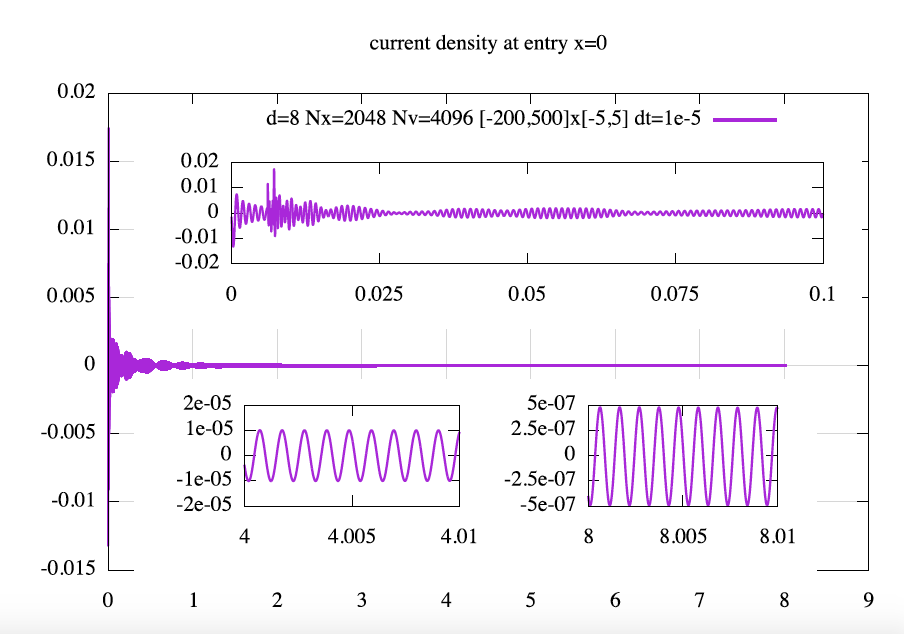} 

\caption{Time evolution of current density $J(x,t)$ at the entry $x=0$ and zoom on time intervals $[0,0.1]$, $[4,4.01]$ and $[8,8.01]$. Parameters given in \eqref{parameters}.}
\label{fig6}
\end{figure}

\textbf{Current density at the entry.} On Figure \ref{fig6}, we represent the  current density at the entry $x=0$ (which is zero at the continuous level). We expect a better behavior on the this boundary than at the wall, where the convergence is more delicate. 
We remark  an oscillatory behavior due to the electron dynamics. After a \emph{violent} behavior at time around $0.01$, where the value can have a peak around $0.02$ (that value is guessed to be linked to the space discretization), it tends to decrease and stabilize, while still oscillating around $0$, with the same frequency.

\begin{figure}
\centering
\includegraphics[width=0.5\linewidth]{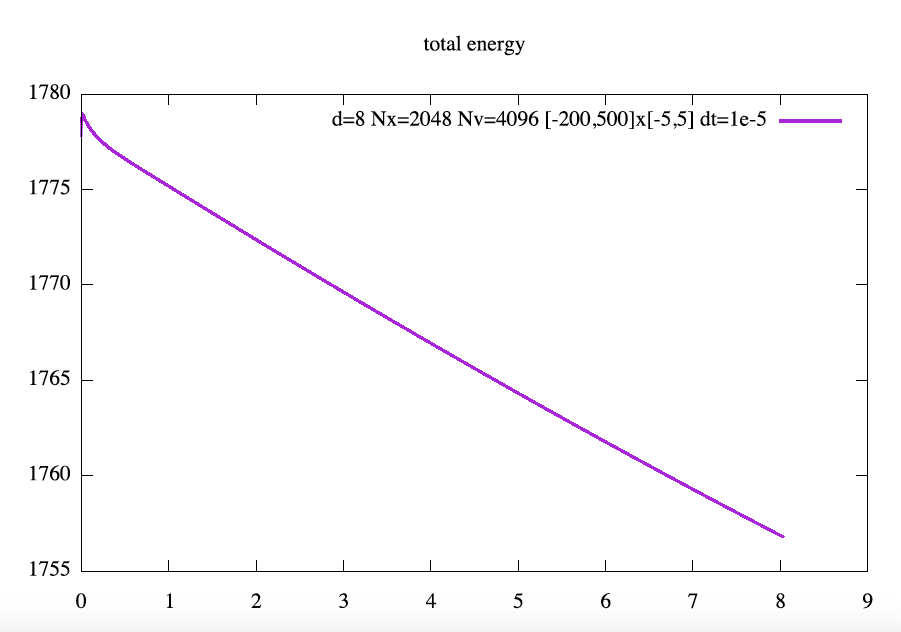} 

\caption{Time evolution of total energy $\mathcal{E}(t)$. Parameters given in \eqref{parameters}.}
\label{fig7}
\end{figure}

\begin{figure}
\begin{tabular}{cc}
\includegraphics[width=0.5\linewidth]{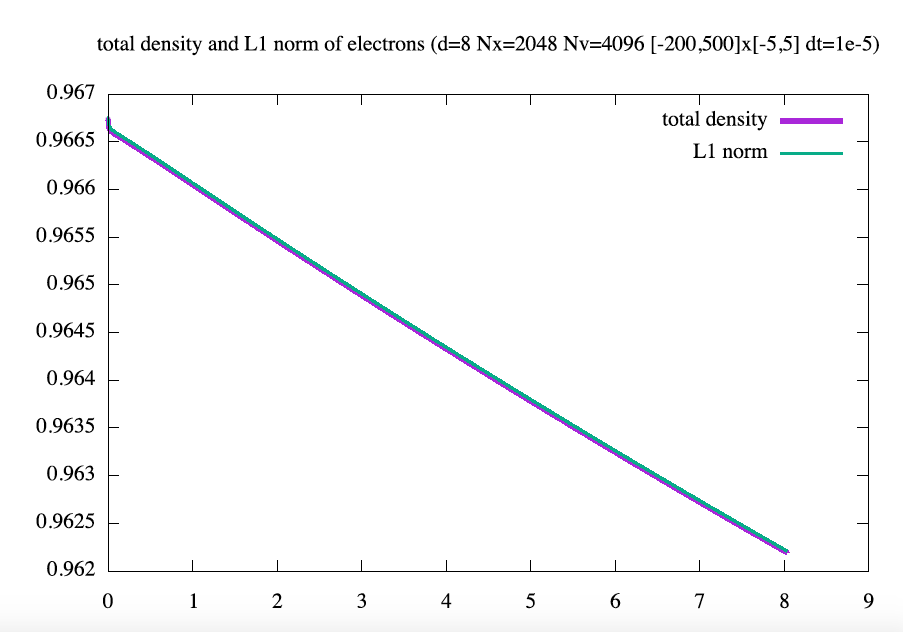}&
\includegraphics[width=0.5\linewidth]{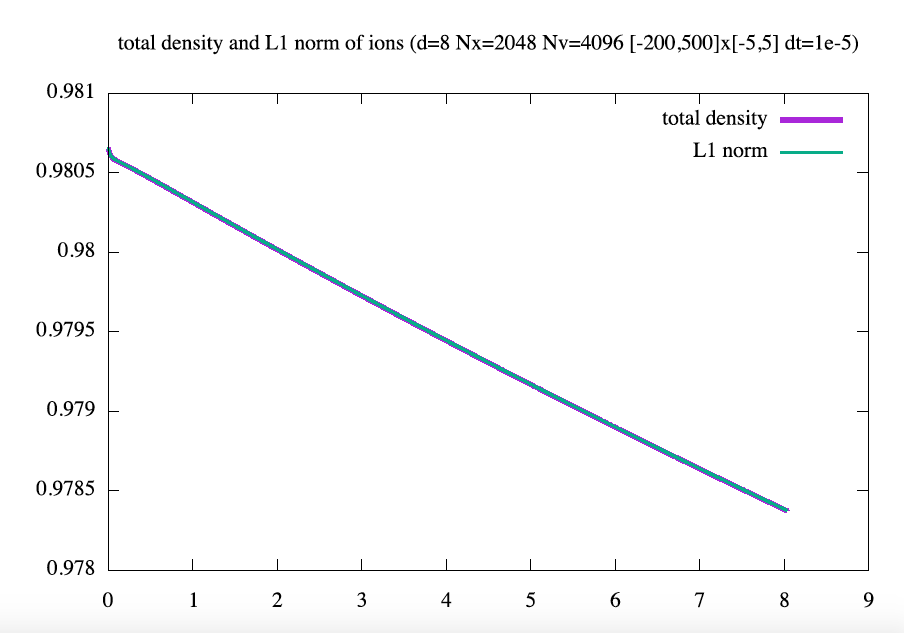}
\end{tabular}
\caption{Time evolution of total density and $L^1$ norm of the distribution function for electrons (left) and ions (right). Parameters given in \eqref{parameters}.}
\label{fig8}
\end{figure}

\begin{figure}
\begin{tabular}{cc}
\includegraphics[width=0.5\linewidth]{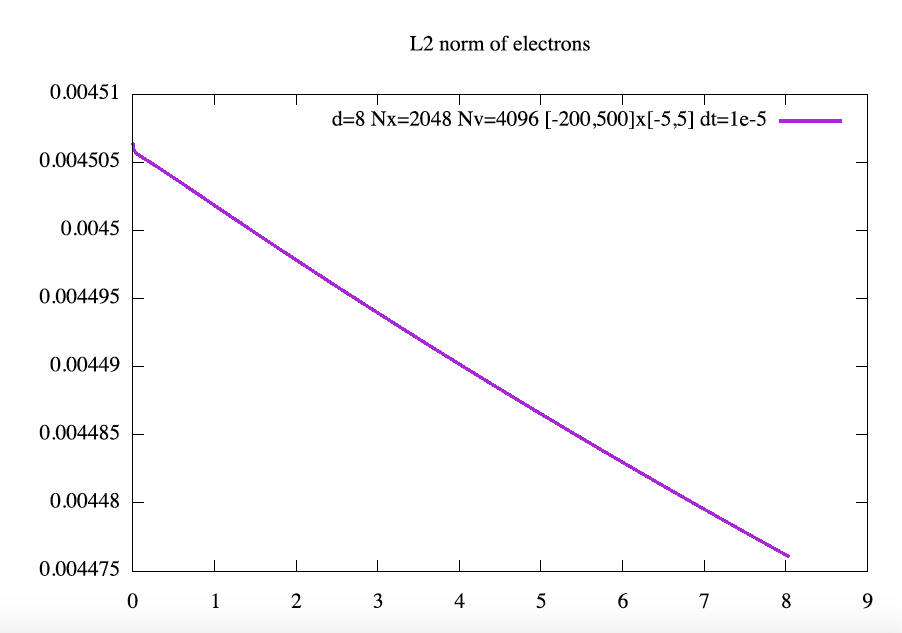}&
\includegraphics[width=0.5\linewidth]{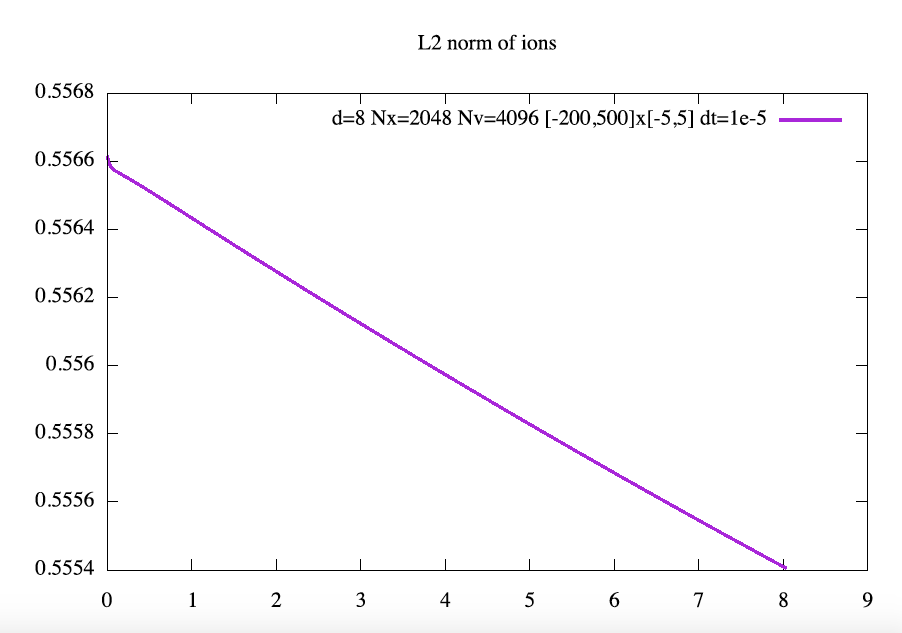}
\end{tabular}
\caption{Time evolution of $L^2$ norm of the distribution function for electrons (left) and ions (right). Parameters given in \eqref{parameters}.}
\label{fig9}
\end{figure}

\textbf{Other diagnostics.}
On Figure \ref{fig7}, we represent the total energy, whose continuous expression is
\begin{equation*}
\mathcal{E}(t) = \frac{1}{2}\int_{v\in \mathbb{R}}\int_0^1 v^{2}\, (f_{s_e}(t,x,v)+f_{s_i}(t,x,v))\, dx dv+\frac{1}{2}\int_0^1 E^2(t,x)\, dx.
\end{equation*}
We see that it is not exactly conserved. Except at the very beginning, it is decreasing and reaches a relative error of about $1.4\%$ at final time.
We then represent the time evolution of total density and $L^1$ norm of the distribution function\footnote{Total density is given by: $\int_{v\in \mathbb{R}}\int_0^1 f_{s}(t,x,v)\, dx dv = \int_0^1 n_{s}(t,x)\, dx$. \\
$L^{1}$ norm of the distribution function is given by: $\int_{v\in \mathbb{R}}\int_0^1 |f_{s}(t,x,v)|\, dx dv$. } for ions and electrons on Figure \ref{fig8} and the $L^2$ norm of the distribution function on Figure \ref{fig9}. We see that their time behavior is similar to the total energy one. Note that the $L^1$ norm and total density are indistinguishable (for both ions and electrons), which shows that the positivity of the distribution functions is rather well preserved.




\subsection{Comparison with other parameters}

\begin{figure}
\begin{tabular}{cc}
\includegraphics[width=0.5\linewidth]{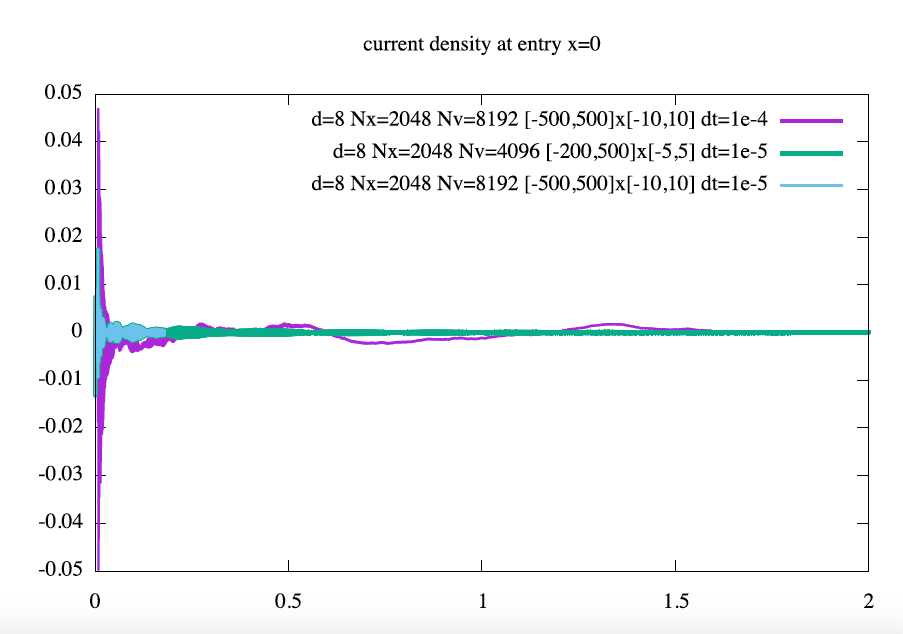}&
\includegraphics[width=0.5\linewidth]{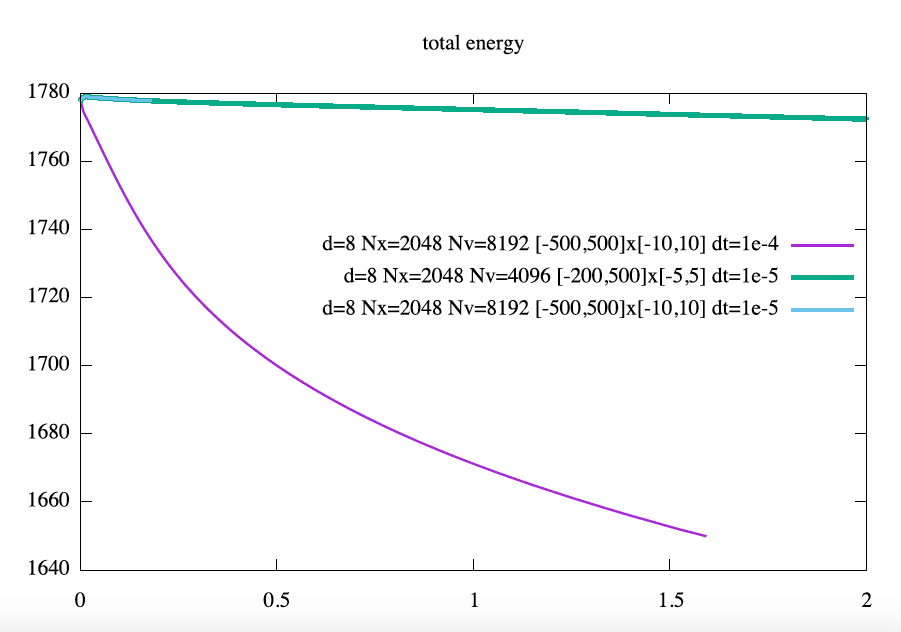}
\end{tabular}
\caption{Comparison of the current density at the entry $x=0$ (left) and total energy (right) changing $\Delta t=10^{-5}$ to  $\Delta t=10^{-4}$.}
\label{fig10}
\end{figure}

On Figure \ref{fig10}, we show how the numerical results change when enlarging the computational domain in velocity (taking velocity domain $[-500,500]$ for electrons and $[-10,10]$ for ions and $N_v=8192$), and then increasing the time step ($\Delta t=10^{-4}$ instead of $\Delta t = 10^{-5}$). We compare the current density at the entry $x=0$ and the total energy. We observe that the change in velocity domain has no influence, while the change in time step $\Delta t$ leads to less accurate results.

\begin{figure}
\begin{tabular}{cc}
\includegraphics[width=0.5\linewidth]{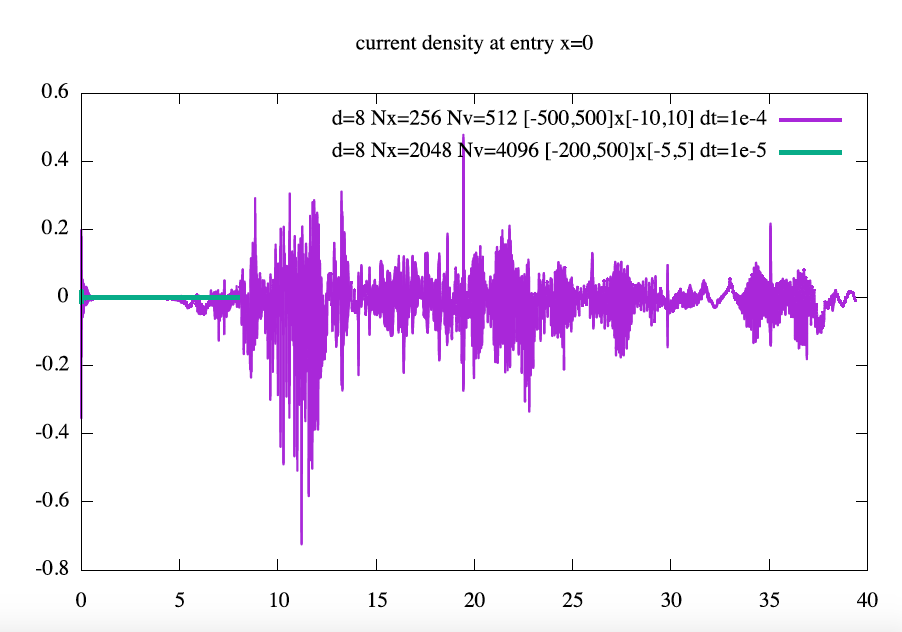}&
\includegraphics[width=0.5\linewidth]{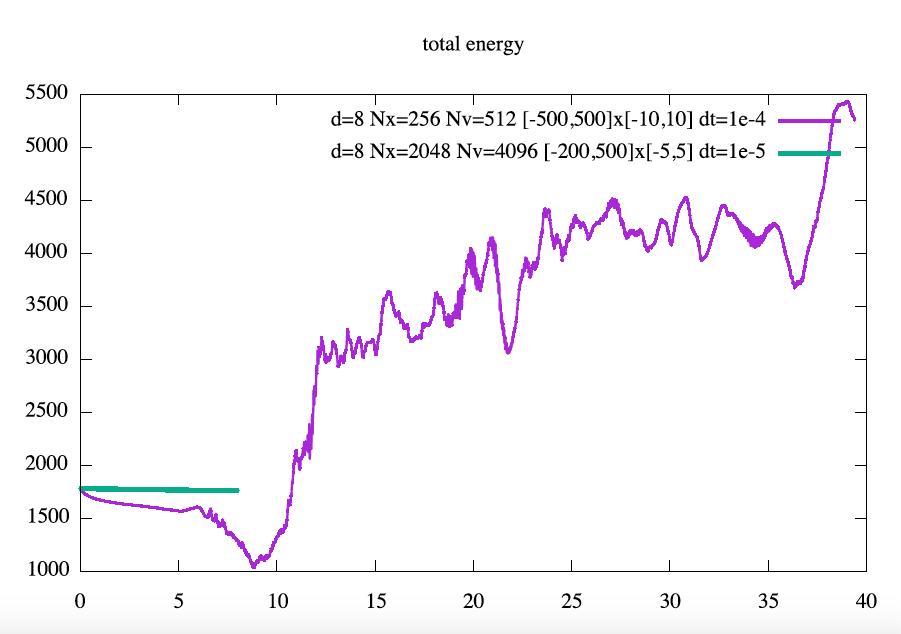}
\end{tabular}
\caption{Comparison of the current density at the entry $x=0$ (left) and total energy (right) changing $\Delta t=10^{-5}$ to  $\Delta t=10^{-4}$ and $N_x$ to $256$ and $N_v$ to $512$.}
\label{fig11}
\end{figure}

On Figure  \ref{fig11}, we present the results when taking less discretization points $N_x=256$, $N_v=512$ and velocity domain $[-500,500]$ for electrons and $[-10,10]$ for ions and $\Delta t=10^{-4}$. This enables to consider longer time simulation, but we 
see that the results are degraded especially in large time.

\begin{figure}
\begin{tabular}{cc}
\includegraphics[width=0.5\linewidth]{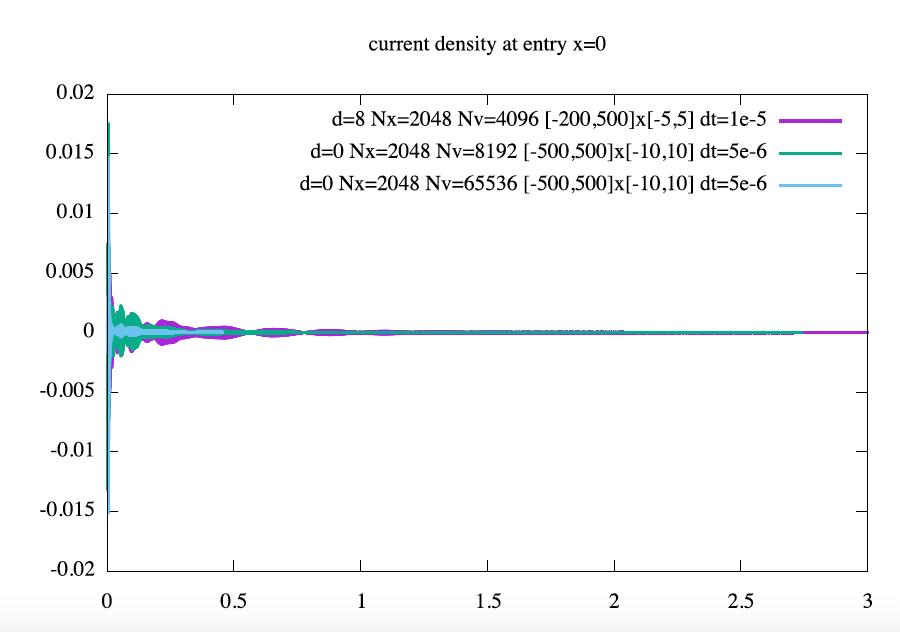}&
\includegraphics[width=0.5\linewidth]{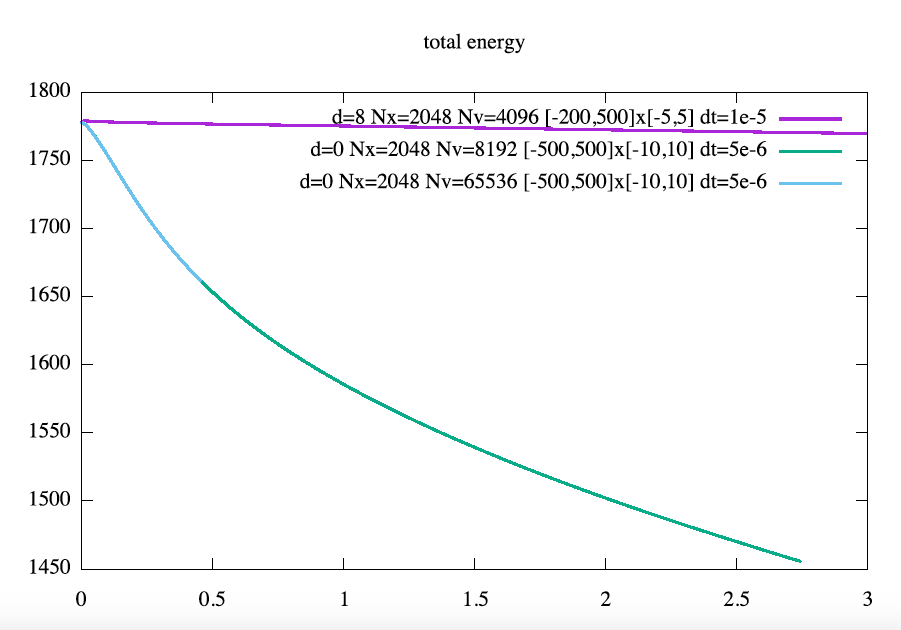}
\end{tabular}
\caption{Comparison of the current density at the entry $x=0$ (left) and the total energy (right) for interpolation degrees $d=8$ and  $d=0$.}
\label{fig12}
\end{figure}

On Figure \ref{fig12}, we consider low order interpolation $d=0$. We consider $\Delta t=5\cdot 10^{-6}$ and $N_v=8192$ or $N_v=65536$. The velocity domain are $[-500,500]$ for electrons and $[-10,10]$ for ions.
We see that the current  density at the entry $x=0$ is improved when using more points in velocity. This is however not true as regards the other diagnostics: for the total energy, the case $d=0$ is really more diffusive and can not compete
with the high order scheme.

\begin{figure}
\begin{tabular}{cc}
\includegraphics[width=0.5\linewidth]{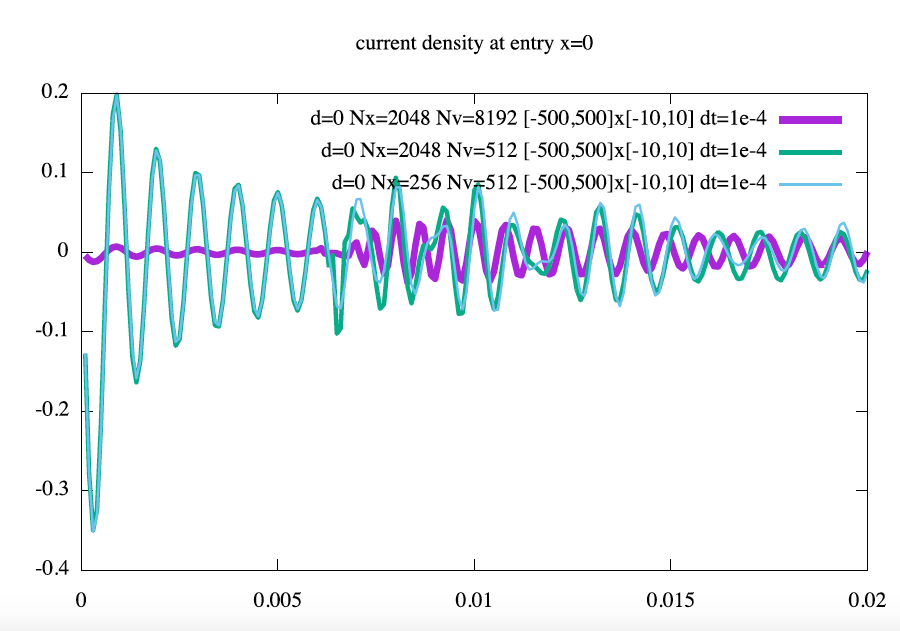}&
\includegraphics[width=0.5\linewidth]{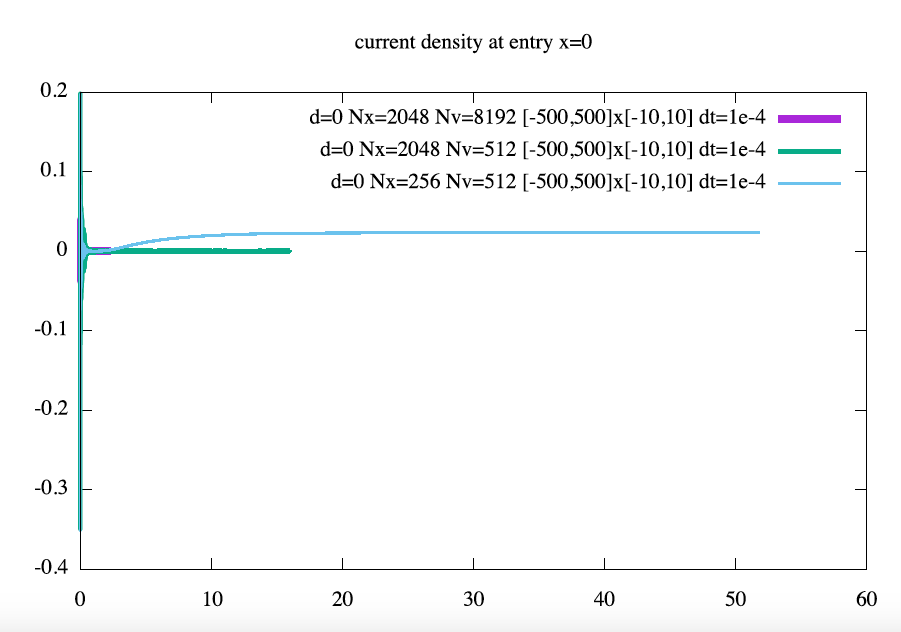}
\end{tabular}
\caption{Current density at the entry $x=0$ for $d=0$, $\Delta t=10^{-4}$ and different space/velocity parameters (left: on short time; right: on long time)}
\label{fig13}
\end{figure}

\begin{figure}
\begin{tabular}{c}
\includegraphics[width=0.5\linewidth]{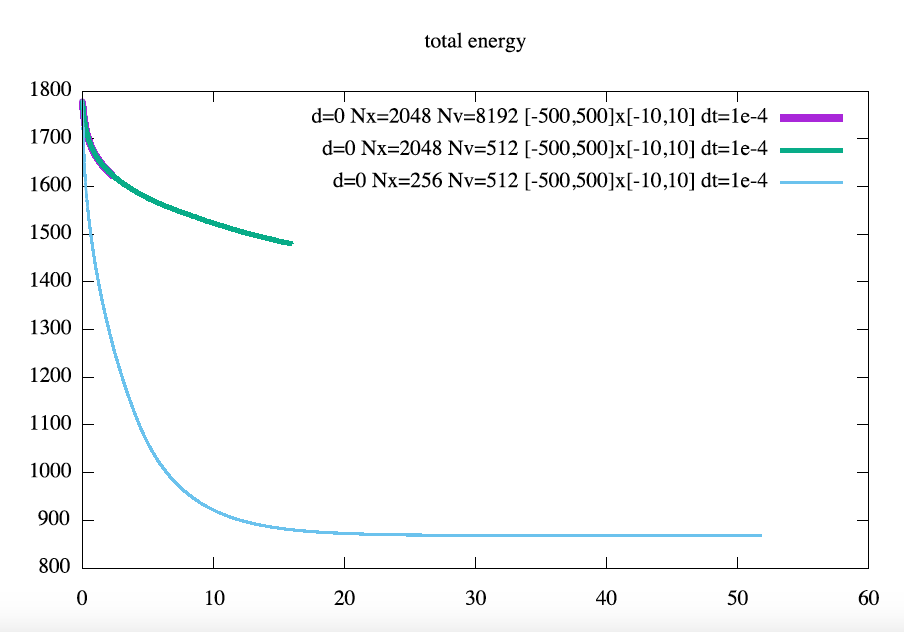}
\end{tabular}
\caption{Total energy for $d=0$, $\Delta t=10^{-4}$ and different space/velocity parameters.}
\label{fig14}
\end{figure}

 On Figures \ref{fig13} and \ref{fig14}, we consider the case $d=0$, with $\Delta t=10^{-4}$ and $N_x=256$ or $N_x=2048$ and $N_v=8192$ or $N_v=512$. Comparing with the previous case with $d=8$ (Figure \ref{fig12}), we see that there is no such oscillatory behavior in large time. However, for $N_x=256$ and $N_v=512$, we see that diffusion is really important and leads to different behavior: the current density at entry tends to a value around $0.1$ instead of $0$.

\section{Conclusion}

We have studied the behavior of the numerical solution of the Vlasov equation, initialized with a sheath equilibrium \cite{2016_Badsi}. 
Thanks to high resolution in velocity, high order interpolation and very small time steps, we are able to recover the equilibrium accurately for relatively long time.
This work is a first step as regards the numerical method. We mention here several directions of research:
\begin{itemize}
\item to reduce the constraint on the time step, asymptotic preserving schemes could be designed, the fictitious boundary values may be improved and we could consider unsplit time integration,
\item to reduce the constraint on space/velocity grid,  adaptive/Discontinuous-Galerkin method or delta-f method could be considered,
\item scheme ensuring a discrete Gauss law could be develop and its impact on the numerical results could be analyzed,
\item we could enhance mixed openmp/mpi parallelization to get full performance on current and future architectures.
\end{itemize}

\section*{Acknowledgment}
This work comes from discussions during the SMAI conference. We thank all the organizers of the conference and also all the speakers of the MULTIKIN mini-symposium.
Thanks to David Coulette, for the discussion and suggestion concerning the boundary conditions.
Thanks to the Selalib team \url{http://selalib.gforge.inria.fr/} and in particular Pierre Navaro.

This work has been carried out within the framework
of the EUROfusion Consortium and has received funding
from the Euratom Research and Training Programme 2014-
2018 under Grant Agreement No. 633053. Computing
facilities were provided by the EUROfusion Marconi supercomputer facility.
The views and opinions expressed herein do not necessarily reflect those of
the European Commission.

\bibliographystyle{plain}
\bibliography{sheath}

\begin{thebibliography}{1}

\bibitem{Afeyan2014}
Bedros Afeyan, Fernando Casas, Nicolas Crouseilles, Adila Dodhy, Erwan Faou,
  Michel Mehrenberger, and Eric Sonnendr{\"u}cker.
\newblock Simulations of kinetic electrostatic electron nonlinear (keen) waves
  with variable velocity resolution grids and high-order time-splitting.
\newblock {\em The European Physical Journal D}, 68(10):295, Oct 2014.

\bibitem{BADSI2017954}
Mehdi Badsi.
\newblock Linear electron stability for a bi-kinetic sheath model.
\newblock {\em Journal of Mathematical Analysis and Applications}, 453(2):954
  -- 972, 2017.

\bibitem{2016_Badsi}
Mehdi Badsi, Martin Campos~Pinto, and Bruno Despr{\'e}s.
\newblock A minimization formulation of a bi-kinetic sheath.
\newblock {\em Kinet. Relat. Models}, 9(4), 2016.

\bibitem{CouletteManfredi2014}
David Coulette and Giovanni Manfredi.
\newblock An {E}ulerian {V}lasov code for plasma-wall interactions.
\newblock {\em Journal of Physics: Conference Series}, 561, 2014.

\bibitem{SLDG-CMV2010}
Nicolas Crouseilles, Michel Mehrenberger, and Francesco Vecil.
\newblock {D}iscontinuous {G}alerkin semi-{L}agrangian method for
  {V}lasov-{P}oisson.
\newblock {\em ESAIM: Proc.}, 32:211--230, 2011.

\bibitem{Sonnen}
Eric Sonnendr\"ucker.
\newblock {\em Numerical methods for the {V}lasov equations.}

\end{thebibliography}
\nocite{*}

\end{document}